\documentclass[a4paper,12pt]{elsarticle}
\usepackage{amssymb}
\usepackage{epsfig}
\usepackage{amsfonts}
\usepackage{amsmath}
\usepackage{euscript}
\usepackage{amscd}
\usepackage{amsthm}
\DeclareMathAlphabet{\mathpzc}{OT1}{pzc}{m}{it}

\newdefinition{defn}[thm]{Definition}
\newdefinition{ex}[thm]{Example}
\newdefinition{rem}[thm]{Remark}
\newdefinition{note}{Note}
\newdefinition{q}{Question}

\newcommand{\comment}[1]{}

\begin{document}
\begin{frontmatter}
\title{A note on factorial $\mathbb{A}^1$-forms with retractions}

\author{Prosenjit Das}
\ead{prosenjit.das@gmail.com}
\address{Stat-Math Unit, Indian Statistical Institute, \\
203, B.T. Road, Kolkata 700 108, India.}

\begin{abstract}
Let $k$ be a field. In this paper we will show that any factorial $\mathbb{A}^1$-form $A$ over any $k$-algebra $R$ is trivial if $A$ has a retraction to $R$.
\smallskip

\end{abstract}
\end{frontmatter}

Throughout this paper, by ``ring'', we shall mean ``commutative ring with unity''.

\medskip

\noindent {\bf Definitions.}
Let $k$ be a field, $R$ a $k$-algebra and $A$ an $R$-algebra. $A$ is said to be an  \textit {$\mathbb{A}^1$-form} over $R$ (with respect to $k$) if there exists an algebraic field extension $k'|_k$ such that $A\otimes _k k' \cong (R \otimes_k k' )[X]$. An $\mathbb{A}^1$-form is called \textit{purely inseparable} (resp. \textit{separable}), if we can take the extension $k'|_k$ to be a purely inseparable (resp. a separable) extension. 

A ring homomorphism $\Phi : A \longrightarrow R$ is said to be a \textit{retraction} if $\Phi$ is $R$-linear. If a retraction $\Phi : A \longrightarrow R$ exists, we say $R$ is a \textit{retract} of $A$.

\medskip

It is well known that any separable $\mathbb{A}^1$-form over any field is trivial. More generally, it has been shown that a separable $\mathbb{A}^1$-form over an arbitrary commutative algebra is trivial (\cite{D_SEP}, Theorem 7), i.e., 

\smallskip

\noindent {\bf Theorem 1.}
\textit{Let $k$ be a field, $L$ a separable field extension of $k$, $R$ a $k$-algebra and $A$ an $R$-algebra such that $A \otimes_k L \cong Sym_{(R \otimes _k L)} (P')$ for a finitely generated rank one projective module $P'$ over $R \otimes _k L$. Then $A \cong Sym_R (P)$ for a finitely generated rank one projective module $P$ over $R$.}

\medskip

If $k$ is not perfect, there exist non-trivial purely inseparable $\mathbb{A}^1$-forms. Asanuma gave a complete structure theorem for purely inseparable $\mathbb{A}^1$-forms over a field $k$ of characteristic $p > 2$ (\cite{A_INSEP}, Theorem 8.1). However, from Asanuma's results, it can be deduced that any factorial $\mathbb {A}^1$-form over a field $k$ with a $k$-rational point is trivial, i.e.,

\smallskip

\noindent {\bf Theorem 2.}
\textit{Let k be a field and $A$ an $\mathbb{A}^1$-form over $k$ such that
\begin{enumerate}
 \item [\rm (1)] $A$ is a UFD.
 \item [\rm (2)] $A$ has a $k$-rational point.
\end{enumerate}
Then $A \cong k[X]$.}

\smallskip

In this paper we shall show that Theorem 2 has a generalization, in the spirit of Theorem 1, to $\mathbb{A}^1$-forms over $k$-UFDs (see Corollary 4 for precise formulation). The generalization turns out to be a special case of our main theorem (Theorem 3) which may be envisaged as an example of ``faithfully flat descent''. 

\comment{\medskip

We register the following property of polynomial algebra which is easy to prove:

\smallskip

\noindent {\bf Lemma.}
\textit{Let $R$ be a ring and $X$ a transcendental element over $R$. Then $ht (XR[X]) = 1$.}

\begin{proof}
Note that if $P$ is a minimal prime ideal of $R[X]$, then $P = P_0 [X]$ for some minimal prime ideal $P_0$ of $R$. Now if $P_0 [X]$ is a minimal prime ideal in a chain of prime ideals determining the height of $XR[X]$, then $ht (XR[X]) = ht (X (R/P_0)[X])$. Since $\underset{n \ge 1} {\overset{} {\cap}} {X^n (R/P_0)[X]} = 0$ and since $R/P_0$ is a domain, $ht (X (R/P_0)[X]) = 1$. Thus $ht (XR[X]) = 1$.
\end{proof}

\smallskip

We now prove the main theorem.}
\smallskip

\noindent {\bf Theorem 3.}
\textit{Let $R$ be a ring and $A$ be an $R$-algebra such that
\begin{enumerate}
 \item [\rm (1)]$A$ is a UFD.
 \item [\rm (2)]There is a retraction $\Phi : A \longrightarrow R$.
 \item [\rm (3)]There exists a faithfully flat ring homomorphism $\eta: R \longrightarrow R'$ such that $A \otimes_R R'\cong R'[X]$.
\end{enumerate}
Then $A \cong R[X]$.}

\begin{proof}
Let $A'=A \otimes_R R'=R'[f]$, $P=Ker \ \Phi$, $P'=PA'(=P \otimes_R R')$,
and let $\Phi'=\Phi \otimes 1$ be the induced retraction from $A'$ to $R'$.
Then $P$ is a prime ideal of $A$ and we have a short exact sequence of $R$-modules 
\[ 0 \longrightarrow P \longrightarrow A \stackrel{\Phi}{\longrightarrow} R \longrightarrow 0 \]
\noindent and hence a short exact sequence of $R'$-modules 
\[0 \longrightarrow P' \longrightarrow A'(=R'[f]) \stackrel{\Phi'}\longrightarrow R' \longrightarrow 0. \]

Then $P' = Ker \ \Phi' = (f - \Phi'(f)) R'[f]$. Replacing $f - \Phi'(f)$ by $f$, we assume $P' = f A' = f R'[f]$. Since $A'$ is faithfully flat over $A$, going-down theorem holds between $A$ and $A'$ (\cite{Mat_RING}, Pg. 68, Theorem 9.5) and also $P' \cap A = P$ (\cite{Mat_RING}, Pg. 49, Theorem 7.5). As $ht (P')=1$, it follows that $ht (P)=1$.
\smallskip

Now, since $A$ is a UFD, there exists $g \in A$ such that $P = gA$. Thus we get $gA' = P' = fA'$. Since $f$ is a non-zero divisor in $A' (= R'[f])$, it follows that $g = \lambda f$ for some unit $\lambda$ in $A'$. Let \[ \lambda = a_0 + a_1 f + a_2 f^2 + \cdots + a_n f^n \] where $a_0$ is a unit in $R'$ and $a_i$ is nilpotent in $R'$ for $1 \le i \le n$. Let $I=(a_1, a_2, \cdots , a_n)R'$. Then $I$ is a nilpotent ideal of $R'$. Let $N$ be the least positive integer such that $I^N=(0)$. Since $g \equiv a_0 f (\hspace{-.25cm}\mod I)$, we have 
\[R'[f]=R'[g] + IR'[f] = \dots = R'[g] + I^N R'[f] = R'[g].\]

Thus we have $R[g] \subseteq A$ and $R[g] \otimes_R R'=R'[g]=A \otimes_R R'$. Since $R'$ is faithfully flat over $R$, it follows that $A=R[g] \cong R[X]$.
\end{proof}

As an immediate consequence of Theorem 3, we get the following corollary:

\smallskip
\noindent {\bf Corollary 4.}
\textit{Let $k$ be a field, $R$ a $k$-algebra and $A$ an $R$-algebra such that
\begin{enumerate}
 \item [\rm (1)]$A$ is a UFD.
 \item [\rm (2)]$R$ is a retract of $A$.
 \item [\rm (3)]$A$ is an $\mathbb{A}^1$-form over $R$.
\end{enumerate}
Then $A \cong R[X]$.}

\medskip

The following two well-known examples (\cite{Unipotent}, Pg. 70--71, Remark 6.6(a), Examples (i) and (ii)) respectively show that in Theorem 2 (and hence in Theorem 3), the hypothesis on the existence of a retraction and the hypothesis ``$A$ is a UFD'' are necessary.

\smallskip

\noindent {\bf Example 1.}
Let $\mathbb{F}_p$ be the prime field of characteristic $p$ and let $k=\mathbb{F}_p(t,u)$ be a purely transcendental extension of $\mathbb{F}_p$ with variables $t$ and $u$. Then $A=k[X,Y]/(Y^p-t-X-uX^p)$ is a factorial non-trivial $\mathbb{A}^1$-form over $k$ which does not have a retraction to $k$.

\smallskip

\noindent {\bf Example 2.}
Let $k$ be a field of characteristic $p \ge 2$ and $A=k [X,Y]/(Y^p-X-aX^p)$ where $a \in k \backslash k^p$. Then $A$ is a non-trivial $\mathbb{A}^1$-form over $k$ with a retraction to $k$. Here $A$ is not a UFD.

\medskip

{\noindent}{\bf Acknowledgement:}

\smallskip
In the earlier version of the article, Corollary 4 was presented as the main theorem. I thank the referee for his observation
that the arguments give the more general statement Theorem 3. 

\smallskip

I thank Amartya K. Dutta for suggesting the problem and for his useful comments and suggestions. I am grateful to T. Asanuma for his lectures on purely inseparable $\mathbb{A}^1$-forms in Indian Statistical Institute, Kolkata. I also thank S. M. Bhatwadekar and N. Onoda for helpful discussions.

\bibliography{reference}
\bibliographystyle{amsalpha}

\end{document}